\documentstyle[11pt,amscd,amssymb]{amsart}
\let\oldsection=\section
\renewcommand{\section}[1]{\vspace{.18in}\par\noindent\addtocounter{section}{1}\setcounter{subsection}{0}{\bf\thesection\hspace{9pt}#1}}
\renewcommand{\subsection}{\par\vspace{.18in}\noindent\addtocounter{subsection}{1}\setcounter{equation}{0}{\bf\thesubsection\hspace{9pt}}}

\theoremstyle{plain}
\newtheorem{thm}{\addtocounter{equation}{1}({\theequation}) Theorem}
\newtheorem{cor}{\addtocounter{equation}{1}({\theequation}) Corollary}
\newtheorem{lem}{\addtocounter{equation}{1}({\theequation}) Lemma}
\newtheorem{prop}{\addtocounter{equation}{1}({\theequation}) Proposition}

\newtheorem{quest}{\addtocounter{equation}{1}({\theequation}) Question}
\theoremstyle{definition}

\numberwithin{equation}{subsection}

\newcommand{\Hom}{\operatorname{Hom}}
\newcommand{\End}{\operatorname{End}}
\newcommand{\ind}{\operatorname{ind}}
\newcommand{\sgn}{\operatorname{sgn}}

\newcommand{\rad}{\operatorname{rad}}

\newcommand{\res}{\operatorname{res}}

\newcommand{\0}{\bar 0}
\newcommand{\1}{\bar 1}
\newcommand{\p}[1]{\overline{#1}}
\newcommand{\Z}{\mathbb{Z}}
\newcommand{\Zt}{{\mathbb{Z}}_{2}}
\newcommand{\Zp}{{\mathbb{Z}}_{p}}
\newcommand{\modp}{\:(\text{mod }p)}

\def\reptypeordinary{1.2.1}
\def\reptypesuper{1.3.1}

\def\equivlemma{2.1.2}
\def\typeMQ{2.2.3}

\def\descenttheorem{3.2.1}
\def\descentforschur{3.3.1}
\def\tensorbyonedim{3.4.1}
\def\ab2trick{3.5.1}
\def\a0pcase{3.6.1}
\def\S11regrep{4.1.1}
\def\answerforS11{4.1.3}

\def\nonprojsigneddistinct{4.2.2}

\def\answerfordlessthan2p{4.2.7}
\def\caseA{4.4.2}
\def\caseB{4.5.2}
\def\caseC{4.6.2}

\def\globaldimension{5.1.1}
\def\signedblocks{5.2.9}

\begin{document}
\title[Representation type of $S(m|n,d)$]
{\bf Representation type of Schur superalgebras}
\author{\sc David J. Hemmer}
\address
{Department of Mathematics\\ University of Toledo \\
2801 W. Bancroft\\Toledo, OH~43606, USA}
\thanks{Research of the first author was supported in part by NSF
grant  DMS-0102019} \email{david.hemmer@@utoledo.edu}
\author{\sc Jonathan Kujawa}
\address
{Department of Mathematics\\ University of Georgia \\
Athens\\ GA~30602, USA}
\thanks{Research of the second author was
supported in part by NSF grant
DMS-0402916}\email{kujawa@@math.uga.edu}
\author{\sc Daniel K. Nakano}
\address
{Department of Mathematics\\ University of Georgia \\
Athens\\ GA~30602, USA}
\thanks{Research of the third author was supported in part by NSF
grant  DMS-0400548} \email{nakano@@math.uga.edu}
\date{October 2004}
\subjclass{Primary 20C30}
\begin{abstract} Let $S(m|n,d)$ be the Schur superalgebra whose supermodules
correspond to the polynomial representations of the supergroup
$GL(m|n)$ of degree $d$. In this paper we determine the
representation type of these algebras (i.e. classify the ones
which are semisimple, have finite, tame and wild representation
type). Moreover, we prove that these algebras are in general not
quasi-hereditary and have infinite global dimension.
\end{abstract}
\maketitle

\maketitle \vskip 1cm

\section{\bf Introduction}

\subsection{} A central problem in the study of finite-dimensional algebras
is to understand the structure of the indecomposable modules. As a
first step one would like to know how many indecomposable modules
a given algebra admits. Any finite-dimensional algebra can be
classified into one of three categories: finite, tame, or wild
representation type. A finite-dimensional algebra $A$ has {\it
finite representation type} if $A$ has finitely many
indecomposable modules up to isomorphism. If $A$ is not of finite
type then $A$ is of {\it infinite representation type}. Algebras
of infinite representation type are either of {\it tame
representation type} or of {\it wild representation type}. For
algebras of tame representation type one has a chance of
classifying all the indecomposable modules up to isomorphism.

There have been much progress in the past ten years in determining
the representation type of important classes of finite-dimensional
algebras. The Schur algebras and the $q$-Schur algebra are
finite-dimensional algebras which arise in the representation
theory of the general linear groups and symmetric groups. A
complete classification of the representation type of these
algebras was given by Doty, Erdmann, Martin and Nakano \cite{erd,
dn, demn, en1}. These algebras are fundamental examples of
quasi-hereditary algebras (or equivalently highest weight
categories) which were introduced by Cline, Parshall and Scott
\cite{cps}. The blocks for the ordinary and parabolic BGG category
${\mathcal O}$ for finite-dimensional complex semisimple Lie
algebras are other important examples of highest weight categories. Results
pertaining to the representation type for these blocks were
obtained independently by Futorny, Nakano and Pollack \cite{fnp},
and Br\"ustle, K\"onig and Mazorchuk \cite{bkm} for the ordinary
category ${\mathcal O}$ and more recently for the parabolic
category ${\mathcal O}$ by Boe and Nakano \cite{bn}. For
quasi-hereditary algebras the projective modules admit filtrations
by certain standard modules. This information allows one to
determine the structure of the projective modules which in turn
can lead to an expression of the algebra via quiver and relations.
In many cases this information can be used to determine the
representation type of the algebra.

Other important classes of finite-dimensional algebras are quasi-Frobenius algebras where all projective modules are injective.
These algebras all have infinite global dimension (except when the
algebra is semisimple). For group algebras, restricted enveloping
algebras, and Frobenius kernels homological information involving
the theory of complexity and support varieties can be used to
deduce information about the representation type of these algebras. For Hecke
algebras of type $A$, the representation type of the blocks was
given by Erdmann and Nakano in \cite{en2}. Ariki and Mathas
\cite{am} classified the representation type for Hecke algebras of
type $B$ using Fock space methods. Ariki \cite{a} recently has
extended this classification to Hecke algebras of classical type.

The main result of this paper is a complete classification of the
representation type for the Schur superalgebras $S(m|n,d)$.  These
algebras are of particular interest because the category of
$S(m|n,d)$-supermodules is equivalent to the category of
polynomial representations of degree $d$ for the supergroup
$GL(m|n).$  Also, as with the classical Schur algebra, there is a
Schur-Weyl duality between $S(m|n,d)$ and the symmetric group
$\Sigma_{d}$ (see \cite{bk}).  We shall show that in most cases
$S(m|n,d)$ is not quasi-hereditary, which is contrary to a recent
conjecture given by Marko and Zubkov. Since the Schur
superalgebras are not quasi-hereditary the aforementioned
filtration techniques cannot be used. In order to obtain
information about the basic algebra for certain Schur superalgebras
we will compute the endomorphism algebras of direct sums of signed
Young modules. This will entail knowing the structures of certain
signed Young modules, which is of independent interest from the viewpoint
of the representation theory of symmetric groups.

\subsection{} We review below the results for the ordinary Schur
algebra $S(m,d)$, although the proof for $S(m|n,d)$ does not rely
on these results. The representation type of $S(m,d)$ was
determined in \cite{erd, dn, demn}\footnote{The statement that $S(2,11)$ for $p=2$ has tame representation type was inadvertently omitted from the original statement of the theorem \cite[Thm. 1.2(A)]{demn}. We thank Alison Parker for pointing this out to us.}.

\begin{thm}
\label{reptypeordinary} Let $S(m,d)$ be the Schur algebra over $k$
where $\operatorname{char }k\geq 0$.
\begin{itemize}
\item[(a)] $S(m,d)$ is semisimple if and only if one of the
following holds
\begin{itemize}
\item[(i)] $\operatorname{char }k=0$; \item[(ii)] $d<p$;
\item[(iii)] $p=2$, $m=2$, $d=3$.
\end{itemize}
\item[(b)] $S(m,d)$ has finite representation type if and only if
one of the following holds
\begin{itemize}
\item[(i)] $p\geq 2$, $m\geq 3$, $d<2p$; \item[(ii)] $p\geq 2$,
$m=2$, $d<p^{2}$; \item[(iii)] $p=2$, $m=2$, $d=5,7$.
\end{itemize}
\item[(c)] $S(m,d)$ has tame representation type if and only if
one of the following holds
\begin{itemize}
\item[(i)] $p=3$, $m=3$, $d=7,8$; \item[(ii)] $p=3$, $m=2$,
$d=9,10,11$; \item[(iii)] $p=2$, $m=2$, $d=4,9,11$.
\end{itemize}
\item[(d)] In all other cases not listed above $S(m,d)$ has wild
representation type.
\end{itemize}
\end{thm}

The classification of the representation type of
the $q$-Schur algebras $S_{q}(m,d)$ can be found in \cite[Thm.
1.4(A)-(C)]{en1}.

\subsection{} Throughout this paper we will assume that $k$ is
an algebraically closed field of characteristic $p\geq 0$. The
following theorem summarizes our results on the Schur
superalgebras.

\begin{thm} Let $S(m|n,d)$ be the Schur superalgebra over $k$ where
$\operatorname{char} k\neq 2$. \label{reptypesuper}
\begin{itemize}
\item[(a)] $S(m|n,d)$ is semisimple if and only if one of the
following holds
\begin{itemize}
\item[(i)] $\operatorname{char} k=0$; \item[(ii)] $d<p$;
\item[(iii)] $m=n=1$ and $p\nmid d$.
\end{itemize}
\item[(b)] $S(m|n,d)$ has finite representation type but is not
semisimple if and only if one of the following holds
\begin{itemize}
\item[(i)] $p\leq d < 2p$; \item[(ii)] $m=n=1$ and $p\mid d$.
\end{itemize}
\item[(c)] In all other cases not listed above $S(m|n,d)$ has wild
representation type.
\end{itemize}
\end{thm}

We should remark that when in characteristic two $S(m|n,d)$ is
equal to the ordinary Schur algebra $S(m+n,d)$, and its
representation type is classified by Theorem~\reptypeordinary.

\section{\bf Comparing Module Categories for Superalgebras}

\subsection{} Recall that $k$ denotes our fixed ground field.  All objects defined over $k$ (algebras, modules, superalgebras, supermodules, etc.) are assumed to be finite-dimensional as $k$-vector spaces.

 A \emph{superspace}
is a $\Zt$-graded $k$-vector space
$V=V_{\0}\oplus V_{\1}.$  Given a homogeneous element of a
superspace, $v \in V,$ we write $\p{v} \in \Zt$ for the degree of
the element.  If $V$ and $W$ are superspaces, then $\Hom_{k}(V,W)$ is naturally a superspace with $\Hom_{k}(V,W)_{\0 }$ (resp.  $\Hom_{k}(V,W)_{\1 }$) consisting of all linear maps which preserve (resp. reverse) grading.

A \emph{superalgebra} is a $\Zt$-graded associative
$k$-algebra $A=A_{\0}\oplus A_{\1}$ which satisfies $A_{r}A_{s} \subseteq A_{r+s}$ for all $r,s \in \Zt.$  An $A$-\emph{supermodule}
is a $\Zt$-graded $A$-module $M=M_{\0}\oplus M_{\1}$ which
satisfies $A_{r}M_{s} \subseteq M_{r+s}$ for all $r,s \in \Zt.$  An $A$-supermodule homomorphism is a linear map $f:M \to N$ which satisfies $f(am)=(-1)^{\p{f}\;\p{a}}af(m)$ for all homogeneous $a \in A$ and all $m \in M.$  As we will do later without comment, the given condition makes sense only for homogeneous elements.  The general case is obtained by linearity.  Note that if $M$ and $N$ are $A$-supermodules, then $\Hom_{A}(M,N)$ is naturally a superspace just as in the previous paragraph.

We can define the \emph{parity change functor} $\Pi$ on the category of $A$-supermodules as follows.  On an object $M$ we have that $\Pi M=M$ as a vector space, but we switch the $\Zt$-grading by setting $(\Pi M)_{r}=M_{r+\1}$ for $r \in \Zt.$  The action of $A$ is defined via $a.m=(-1)^{\p{a}}am$ for homogeneous $a \in A$ and $m \in M.$  On a morphism $f,$ $\Pi f$ is the same linear map as $f.$  We refer the reader to \cite[Sect. 2]{bkl} for further discussion of general results on superalgebras, supermodules, etc.

Throughout the remainder of this subsection we assume
$\operatorname{char }k \neq 2.$  Given a superalgebra $A$ there are three module categories which
are natural to consider:
\begin{itemize}
\item[(i)] the category $A$-smod of all
$A$-supermodules and \emph{all} (not necessarily graded)
$A$-supermodule homomorphisms;
\item[(ii)] the underlying even category of $A$-smod
where we take the objects of $A$-smod but only the \emph{even}
(i.e. grading preserving) homomorphisms;
\item[(iii)] the category $A$-mod of all $A$-modules.
\end{itemize}
Observe that an $A$-supermodule is indecomposable or irreducible regardless of whether we view it as an object of the category $A$-smod or its underlying even category.  If $M$ and $N$ are indecomposable $A$-supermodules and $f: M \to N$ is an $A$-supermodule isomorphism, write $f=f_{\0}+f_{\1}$ with $f_{r} \in \Hom_{A}(M,N)_{r}$ for $r \in \Zt.$  Note that by considering the $\Zt$-grading it is straightforward to verify that $f_{\0}-f_{\1}:M \to N$ is also an isomorphism.  Using this one can show that the map $M \oplus \Pi M \to N \oplus \Pi N$ given by
\[
(m,n) \mapsto (f_{\0}(m)+f_{\1}(n),f_{\1}(m)+f_{\0}(n))
\] is in fact an \emph{even} isomorphism of $A$-supermodules.  Consequently, by the graded version of the Krull-Schmidt theorem we see that $M$ is isomorphic to either $N$ or $\Pi N$ by an even isomorphism.  Taken together the above discussion proves that the category $A$-smod is semisimple, finite, tame, or wild if and only if its underlying even category is semisimple, finite, tame, or wild, respectively.  We use this observation without comment in what follows.
The remainder of this section is
devoted to studying the relationship between the representation
type of these categories and the category $A$-mod.

Given a superspace $V,$ let $\delta: V \to V$ be the linear involution
defined by $\delta(v)=(-1)^{\p{v}}v,$ where $v \in V$ is
homogeneous. In particular, the involution $\delta: A \to A$ is an
algebra homomorphism and defines an action of $\Zt$ on $A.$
Consequently, we can define the smash product of $k\Zt$ and $A$ as
follows.  Set
\begin{equation*}
B= k\Zt \otimes_{k} A
\end{equation*} as a vector space.  Identify $k\Zt$ with $k[x]/(x^{2}-1).$  Then the product in $B$ is given by
\begin{equation}
(x^{r} \otimes a)(x^{s} \otimes b)= x^{r+s} \otimes
\delta^{s}(a)b,
\end{equation}
where $r, s \in\{0,1 \}$ and $a,b \in A.$
We identify $k\Zt$ and $A$ as subalgebras via
$k\Zt \otimes 1 $ and $1 \otimes A,$ respectively.
Note that since we assume $\operatorname{char }k\neq 2,$ $x$ acts
semisimply on any $k\Zt$-module with eigenvalues $\pm 1$.  The following result is proven in \cite[Thm. 1.4b]{Mui}.

\begin{prop}\label{equivlemma}  The category of $B$-modules is isomorphic to the
underlying even category of $A$-supermodules.
\end{prop}

\begin{pf}
 Let $M$ be an $A$-supermodule.  We define an action of
$k\Zt$ on $M$ via $xm=\delta(m).$  It is straightforward to verify
that this makes $M$ into a $B$-module. Conversely, if $M$ is a
$B$-module, then we can view it as an $A$-module by restriction.
To obtain a $\Zt$-grading we set $M_{\0}$ (resp. $M_{\1}$) to be the
$1$ (resp. $-1$) eigenspace of $x$ acting on $M$.  One can verify that this
$\Zt$-grading makes $M$ into an $A$-supermodule.  It is
straightforward to verify that a morphism in one category
restricts or extends to a morphism in the other category.
\end{pf}

\subsection{} Note that there is an involution $\zeta: B \to B$ given by $x^{r}\otimes a \mapsto (-1)^{r}x^{r} \otimes a$ for $r \in \{0,1 \}$ and $a \in A.$  Given a $B$-module $M$ we can twist it by $\zeta,$ which we denote by $M^{\zeta}.$  Namely, $M^{\zeta}=M$ as a vector space but with the action of $B$ given by $b.m=\zeta(b)m$ for all $b \in B$ and $m\in M.$  Note that $M^{\zeta}$ is indecomposable if and only if $M$ is indecomposable.

Since $A$ is a subalgebra of $B,$ there is a restriction functor
$\res_{A}^{B}:B\text{-mod}\to A\text{-mod}$ and an induction
functor $\ind_{A}^{B}:A\text{-mod}\to B\text{-mod} $ given by
\[
\ind_{A}^{B}N=B \otimes_{A}N
\] for any $A$-module $N.$  Note that
\begin{equation}\label{E:inddecomp}
\ind_{A}^{B}N \cong 1 \otimes N \oplus x \otimes N
\end{equation}
as $A$-modules.  If $M$ is a $B$-module,  then we have
\begin{equation}\label{E:indiso}
M \oplus M^{\zeta} \cong \ind_{A}^{B}\res_{A}^{B}M.
\end{equation}
Since $x$ acts semisimply on $M$ and $M^{\zeta},$ it suffices to define the
isomorphism on eigenvectors $m \in M,$ $m' \in M^{\zeta}$ of eigenvalue
$\varepsilon, \varepsilon' \in \{\pm 1 \},$ respectively.  Then
the map $M \oplus M^{\zeta} \to \ind_{A}^{B}\res_{A}^{B}M$ is given by
\[
(m,m') \mapsto (1+\varepsilon x) \otimes m +  (1-\varepsilon' x)\otimes m'.
\]  To verify that the map is $B$-linear, one uses the observation that if $a \in A$ is an eigenvector for the action of $x$ with eigenvalue $\alpha$ and $m \in M$ is an eigenvector for the action of $x$ with eigenvalue $\beta,$ then $am$ is an eigenvector with eigenvalue $\alpha\beta.$

\begin{prop}\label{typeMQ}  Let $M$ be an indecomposable $B$-module.  Then $M$ is either indecomposable as an $A$-module or the direct sum of two indecomposable $A$-modules.   Furthermore, every indecomposable $A$-module $N$ is a direct summand of $\ind_{A}^{B}N$ as $A$-modules.

On the other hand, if $M$ is an indecomposable $B$-module which contains $N$ as a direct summand, then $M$ is unique up to isomorphism and twisting by $\zeta.$

\end{prop}
\begin{pf} Let $M$ be an indecomposable $B$-module and let
\begin{equation}\label{E:directsumdecomp}
M= M_{1} \oplus \dotsb \oplus M_{k}
\end{equation}
be a decomposition of $M$ into indecomposable $A$-modules. By
using \eqref{E:indiso} and applying $\ind_{A}^{B} $ one has
\[
M \oplus M^{\zeta} \cong \ind_{A}^{B}M_{1} \oplus \dotsb \oplus
\ind_{A}^{B}M_{k}.
\]
Since the left hand side is a decomposition as indecomposable
$B$-modules, by the Krull-Schmidt theorem it follows that
\eqref{E:directsumdecomp} cannot have more than two direct
summands.

Now say $N$ is an indecomposable $A$-module.   By
\eqref{E:inddecomp} we see that $N$ is a direct summand of $\ind_{A}^{B} N.$   On the other hand, if $M$ is an indecomposable $B$-module which contains $N$ as a direct summand, then $\ind_{A}^{B}M \cong M \oplus M^{\zeta}$ contains $\ind_{A}^{B}N.$  If $\ind_{A}^{B} N$ is indecomposable, then either $M \cong \ind_{A}^{B}N$ or $M \cong (\ind_{A}^{B}N)^{\zeta}.$  If $\ind_{A}^{B} N$ is decomposable, then $\ind_{A}^{B}N\cong M \oplus M^{\zeta}.$  This proves the desired result.

\end{pf}

As is well known, if $M$ is an irreducible $A$-supermodule then as
an $A$-module $M$ is either irreducible or the direct sum of two
irreducibles.  Combining the previous result with
Proposition~\typeMQ\ shows this behaviour extends to
indecomposables.

\subsection{}  We can now prove the main result of this section: the representation type of $A$ as a superalgebra coincides with that of $A$ as an algebra.

\begin{thm}\label{samereptype}  The category $A$-smod is semisimple, finite, tame, or wild if and only if the same is true for $A$-mod.
\end{thm}

\begin{pf} By Proposition~\equivlemma\ it suffices to consider the category
$B$-mod rather than $A$-smod. \vskip .25cm \noindent {\it Semisimple
type:}  That $A$-smod is semisimple if and only if $B$-mod is semisimple is already known.  See, for example, \cite[Sect. 2]{bkl} or \cite[Thm. 1.4c]{Mui}.
\vskip .25cm \noindent {\it Finite
type:}  It is clear from Proposition~\typeMQ\ that $A$-mod has
finitely many indecomposables if and only if $B$-mod does. \vskip
.25cm \noindent {\it Tame type:}  Suppose $A$-mod has tame
representation type.  Let $t$ be an indeterminant and $d$ a
positive integer. Then there are indecomposable $A$-modules
$D_{d,1}, \dotsc , D_{d,r_{d}}$ and indecomposable
$A-k[t]$-bimodules $E_{d,1}, \dotsc , E_{d,s_{d}}$ which are free
as right $k[t]$-modules so that any indecomposable $A$-module of
dimension $d$ is either isomorphic to $D_{d,i}$ for some $i,$ or
is isomorphic to $E_{d,i} \otimes_{k[t]}L$ for some $i$ and some
irreducible $k[t]$-module $L.$  Observe that the proof
Proposition~\typeMQ\ applies verbatim to indecomposable
$A-k[t]$-bimodules. Consequently, we can choose indecomposable
$B$-modules $\widetilde{D}_{d,1}, \dotsc ,
\widetilde{D}_{d,r_{d}}$ and indecomposable $B-k[t]$-bimodules
$\widetilde{E}_{d,1}, \dotsc , \widetilde{E}_{d,s_{d}}$ (which are
free as $k[t]$-supermodules) where $\widetilde{D}_{d,i}$ contains
$D_{d,i}$ as a direct summand as $A$-modules and
$\widetilde{E}_{d,i}$ contains $E_{d,i}$ as a direct summand as
$A-k[t]$-bimodules.

Let $M$ be an indecomposable $B$-module and $N$ an indecomposable
direct summand of $M$ as an $A$-module.  If $N \cong D_{d,i}$ for
some $d$ and $i,$ then by the uniqueness statement in
Proposition~\typeMQ\ it follows that $M \cong \widetilde{D}_{d,i}$
or $M \cong (\widetilde{D}_{d,i})^{\zeta}.$   On the other hand,
say $N \cong E_{d,i} \otimes_{k[t]} L$ for some $d,i,$ and
irreducible $k[t]$-module $L.$  Since $E_{d,i}$ is a direct
summand of $\widetilde{E}_{d,i},$ we have that $N \cong E_{d,i}
\otimes_{k[t]} L$ is a direct summand of
$\widetilde{E}_{d,i}\otimes_{k[t]} L$.  However, by \cite[Lemma
2.2]{demn} we know $\widetilde{E}_{d,i}\otimes_{k[t]} L$ is
indecomposable.  By uniqueness again we have $M \cong
\widetilde{E}_{d,i}\otimes_{k[t]} L$ or  $M \cong
(\widetilde{E}_{d,i}\otimes_{k[t]}
L)^{\zeta}=(\widetilde{E}_{d,i})^{\zeta}\otimes_{k[t]} L.$
Therefore, if $A$-mod is tame, then $B$-mod is tame. \vskip .25cm
\noindent {\it Wild type:} Suppose that  $A$-mod has wild
representation type. By \cite[Theorem 2.2]{demn} and
Proposition~\typeMQ\ it follows that $B$-mod has wild
representation type.

\end{pf}

\section{\bf Results on the Schur Superalgebra}

\subsection{\bf The Schur Superalgebra:} Fix nonnegative integers
$m,n,$ and $d,$ and fix a superspace $V$ with fixed basis
$v_{1}, \dotsc , v_{m+n}$ and $\Zt$-grading given
by setting $\p{v}_{i}=\0 $ for $i=1, \dotsc , m$ and $\p{v}_{i}=\1
$ for $i=m+1, \dotsc , m+n.$  The tensor space $V^{\otimes d}$ is
then naturally $\Zt$-graded by setting the degree of $v_{i_{1}}
\otimes \dotsb \otimes v_{i_{d}}$ to be $\p{v}_{i_{1}} + \dotsb +
\p{v}_{i_{d}}$ for all $1 \leq i_{1},\dotsc ,i_{d} \leq m+n. $
The symmetric group $\Sigma_{d}$ acts on $V^{\otimes d}$ on the
right via
\begin{equation*}
v_{i_{1}} \otimes \dotsb \otimes v_{i_{j}}\otimes v_{i_{j+1}}\otimes \dotsb  \otimes v_{i_{d}} (j\:\:j+1) :=
(-1)^{\bar{v}_{i_{j}}\bar{v}_{i_{j+1}}} v_{i_{1}}\otimes \dotsb
\otimes v_{i_{j+1}} \otimes v_{i_{j}} \otimes \dotsb \otimes
v_{i_{d}}
\end{equation*}
for $1 \leq i_1,\dots,i_d \leq m+n$ and each $1 \leq j < d$.

The Schur superalgebra is defined to be
\[
S(m|n,d)=\End_{\Sigma_{d}}( V^{\otimes d}),
\] with the $\Zt$-grading given by setting $S(m|n,d)_{\0}$ (resp. $S(m|n,d)_{\1}$) to be the set of all degree preserving (resp. degree reversing) maps.  Observe that when $p=2$ the above action coincides with the usual permutation action and so $S(m|n,d)=S(m+n,d).$

We also note that $S(m|n,d) \cong S(n|m,d)$ as superalgebras.  Namely, let $t: \{1, \dotsc , m+n \} \to \{1, \dotsc, m+n\}$ be given by $t(i)=i+n$ if $i=1, \dotsc , m$ and $t(i)=i-m$ if $i=m+1, \dotsc , m+n.$  Fix a superspace $\widetilde{V}$ with fixed basis
$\tilde{v}_{1}, \dotsc , \tilde{v}_{m+n}$ and $\Zt$-grading given
by setting $\p{\tilde{v}}_{i}=\0 $ for $i=1, \dotsc , n$ and $\p{\tilde{v}}_{i}=\1
$ for $i=n+1, \dotsc , m+n.$  Define a linear map $\tau: V^{\otimes d} \otimes \sgn \to \widetilde{V}^{\otimes d}$ via
\begin{equation*}
\tau(v_{i_{1}} \otimes \dotsb \otimes v_{i_{d}}\otimes 1)=(-1)^{(d-1)\p{v}_{i_{1}}+(d-2)\p{v}_{i_{2}}+\dotsb + \p{v}_{i_{d-1}}} \tilde{v}_{t(i_{1})} \otimes \dotsb \otimes \tilde{v}_{t(i_{d})},
\end{equation*}
where $v_{i_{1}}, \dotsc , v_{i_{d}}$ are elements of our fixed homogeneous basis for $V.$  It is straightforward to verify on simple transpositions that in fact $\tau$ is a $k\Sigma_{d}$-module isomorphism (c.f. the proof of \cite[Theorem 6.4]{bk}).  Consequently we have
\[
S(m|n,d)=\End_{\Sigma_{d}}(V^{\otimes d}) \cong \End_{\Sigma_{d}}(V^{\otimes d}\otimes \sgn ) \cong \End_{\Sigma_{d}}(\widetilde{V}^{\otimes d})=S(n|m,d),
\] and the isomorphisms preserve the $\Zt$-grading.

We now introduce certain modules which will be crucial in later
calculations.  Given a tuple of integers $\gamma=(\gamma_{1}, \dotsc , \gamma_{t})$ we denote the integer $\gamma_{1}+\dotsb +\gamma_{t}$ by $|\gamma|.$  Let
\[
\Lambda(m|n,d)=\left\{ (\lambda_{1}, \dotsc , \lambda_{m} | \mu_1,
\ldots, \mu_n) \in {\Z}^{m+n} \mid \lambda_{i}, \mu_i \geq 0 \text{ for all } i
\text{ and } |\lambda|+|\mu|=d\right\}.
\]
We write $(\lambda | \mu)$ for an element of $\Lambda(m|n,d)$.

In general, given $\gamma=(\gamma_{1}, \dotsc , \gamma_{t}),$ a sequence of nonnegative integers summing to $r,$ let
\[
\Sigma_{\gamma} = \Sigma_{\gamma_{1}} \times
\Sigma_{\gamma_{2}} \times \dotsb \times \Sigma_{\gamma_{t}}
\] viewed as a subgroup of $\Sigma_{r}$ in the natural way by $\Sigma_{\gamma_{1}}$
acting as permutations of the first $\gamma_{1}$ letters,
$\Sigma_{\gamma_{2}}$ acting as permutations of the next
$\gamma_{2}$ letters, and so on.

Given $(\lambda | \mu) \in \Lambda(m|n,d)$ we define the
\emph{signed permutation module} $M^{(\lambda | \mu)}$ as
\[
M^{(\lambda | \mu)}= \ind_{\Sigma_{\lambda}\times
\Sigma_\mu}^{\Sigma_{d}} k \boxtimes \sgn,
\] where $k$ denotes the trivial module for
$\Sigma_{\lambda_{1}} \times \dotsb \times \Sigma_{\lambda_{m}}$
and $\sgn$ denotes the one dimensional sign representation for
$\Sigma_{\mu_1} \times \dotsb \times \Sigma_{\mu_{n}}.$

By definition a $\Sigma_{d}$-module is a \emph{signed Young
module} if for some nonnegative integers $m$ and $n$ and $(\lambda | \mu) \in \Lambda(m|n,d)$ it is a direct summand of $M^{(\lambda | \mu)}.$  Signed permutation modules and signed
Young modules were first considered by Donkin \cite{Do}.  Finally we observe the well known fact that as
$k\Sigma_d$ modules with the action from above we have:

\begin{equation}
\label{tensorspacedecomp}V^{\otimes d} \cong \bigoplus_{(\lambda
|\mu) \in \Lambda(m|n, d)} M^{(\lambda | \mu)}.
\end{equation}

\subsection{} We begin by collecting a few preliminary results.

\begin{prop}
    \label{descenttheorem}
Let $S$ be a finite-dimensional algebra containing an idempotent
$e$, so $eSe$ is a subalgebra of $S$.
\begin{itemize}
\item[(a)] Every indecomposable
$eSe$-module is a summand of the restriction to $eSe$ of some
indecomposable $S$-module.
\item[(b)] If $eSe$ has infinite representation type then so does $S$.
\item[(c)] If $eSe$ has wild representation type then so does $S$.
\end{itemize}
\end{prop}
\begin{pf}
(a) Let $N$ be an indecomposable $eSe$-module. Let $U=S \otimes_{eSe}
N$. Then $U$ is a nonzero $S$-module and:
$$U \cong eU \oplus (1-e)U$$ as $eSe$ modules. But $eU \cong N$.
Now writing $U$ as a direct sum of indecomposable $S$ modules and
an easy argument with the Krull-Schmidt theorem gives the result.
(b) This is clear from part (a).  If $S$ had only finitely many
indecomposable modules then part (a) would imply $eSe$ does as
well. (c) This statement follows from the argument of Theorem 2.2
in \cite{demn}. If the statement of (a) holds for an arbitrary
algebra $A$ and subalgebra $B$, the authors show that if $A$ is
tame then so is $B$.
\end{pf}

\subsection{} One can now apply the previous result to specific situation of $S(m|n,d)$:
This allows us to descend from the superalgebra $S(m'|n', d)$ to
$S(m|n,d)$ for $m \leq m'$, $n \leq n'$.

\begin{cor} \label{descentforschur} Let $m \leq m'$ and $n \leq n'$. Then:
\begin{itemize}
\item[(a)]  If $S(m|n, d)$ has infinite
representation type then so does $S(m'|n', d)$.
\item[(b)] If $S(m|n, d)$ has wild representation type then so does
$S(m'|n', d)$.
\end{itemize}
\end{cor}
\begin{pf} Recall that $\Lambda(m|n,d)$ is the set of pairs
$(\lambda | \mu)$ where $\lambda$ and $\mu$ are compositions with
$|\lambda| + |\mu|=d$ and $\lambda$ has at most $m$ parts, $\mu$
has at most $n$ parts. Thus $\Lambda(m|n,d)$ is naturally a subset
of $\Lambda(m'|n',d)$. Recall that:

$$S(m'|n',d) \cong \End_{k\Sigma_d}\left(\bigoplus_{(\lambda|\mu) \in
\Lambda(m'|n',d)}M^{(\lambda| \mu)}\right).$$

So if we let $e$ denote the projection onto the direct sum of all
$M^{(\lambda| \mu)}$ with $(\lambda| \mu) \in \Lambda(m|n,d)$ then
it is clear that $eS(m'|n', d)e \cong S(m|n,d)$, so we can apply
Proposition ~\descenttheorem\ to prove the result.
\end{pf}

\subsection{}The classical (and infinitesimal) Schur algebra
$S(n,n)$ has a natural one-dimensional representation
corresponding to the determinant for $GL(n)$. This allows one to
embed the module category for $S(n,d)$ into that of $S(n,d+n)$,
which can be used to study representation type (e.g. \cite[Thm.
2.2]{DNPfinitetype}). Unfortunately the analogue of the determinant representation
for the supergroup $GL(m|n),$ the one-dimensional Berezinian representation, is not polynomial and so it does not
give a natural one-dimensional $S(m|n,m+n)$-module.

We instead use
the existence of a one-dimensional representation
for $S(2|1,p)$.  Namely, let $L(0,0|p)$ denote one-dimensional supermodule which is the Frobenius twist of the one-dimensional $GL(2) \times GL(1)$-module $L(0,0) \boxtimes L(1)$ (see \cite[Remark 4.6]{bk}).
Since the polynomial representations of $GL(m|n)$ correspond with the
representations of the supermonoid $M(m|n)$ one can
follow the arguments given in \cite[Prop. 2.1, Prop. 2.2]{DNPfinitetype} by
replacing the determinant representation with $L(0,0|p)$ to deduce
the following result.

\begin{prop}
\label{tensorbyonedim} If $S(2|1, d)$ has wild representation type
then so does $S(2|1,d+cp)$ for any $c \in \mathbf{N}$.
\end{prop}

\subsection{} We will assume familiarity with basic modular representation
theory of the symmetric group. In the next section we will prove
the main result of the paper. Surprisingly, we will use virtually no
representation theory of the superalgebra. However,
we will need to know that Young modules
$Y^{(a, b, 1^c)}$ are actually signed Young modules when $c<p$ and
$m=2, n=1$.

\begin{prop}
\label{ab2trick}Let $1 \leq c <p$. Then the Young module
$Y^{(a,b,1^c)}$ is a summand of the signed permutation module
$M^{(a,b|c)}$.
\end{prop}
\begin{pf}
Since $c<p$ we have

\begin{equation}
\label{regrepc} M^{(1^c)} \cong k\Sigma_c \cong
\ind_{\Sigma_{1^c}}^{\Sigma_c}k \cong \bigoplus_{\lambda \vdash
c}(S^\lambda)^{\dim S^\lambda}.
\end{equation}
But the Specht module $S^\lambda$ has dimension greater than one
except when $\lambda = (1^c)$ or $\lambda=c$. Thus using
(\ref{regrepc}) we obtain:
\begin{equation}
\label{bigperm}
 M^{(a,b,1^c)} \cong M^{(a,b,c)} \oplus M^{(a,b|c)}
\bigoplus_{\lambda \vdash c, \,\,\lambda \neq(1^c),\,\, (c)}
\ind_{\Sigma_a \times \Sigma_b \times \Sigma_c}^{\Sigma_d}(k
\boxtimes k \boxtimes S^\lambda)^{d_\lambda}
\end{equation}
 with all the
$d_\lambda
>1$.
Since $Y^{(a,b,1^c)}$ occurs exactly once as a direct summand of
$M^{(a,b,1^c)}$ and
$\Hom_{\Sigma_d}(Y^{(a,b,1^c)},M^{(a,b,c)})=0,$ the result follows
from (\ref{bigperm}).

\end{pf}

\subsection{} The significance of Proposition~\ab2trick\ is that it gives us a Young
module $Y^\lambda$ as a summand of the tensor space $V^{\otimes d}$, even though $\lambda$ has more than three parts and
$m+n=3$. We need information about one other special signed Young
module for the $p=3$ case.

\begin{prop}
\label{a0pcase} Suppose $p$ does not divide $a$. Then:
$$M^{(a,0|p)} \cong S^{(a+1,1^{p-1})} \oplus S^{(a,1^{p})}.$$
The two Specht modules in the decomposition are both irreducible
and:
$$S^{(a,1^{p})} \cong D^{(a,2,1^{p-2})}.$$
\end{prop}
\begin{pf}
The Specht module decomposition follows from Young's rule (see
\cite[p.51]{J1}). The irreducibility of the ``hook" Specht modules
when $p \nmid d$ is well known work of Peel \cite{P}. Finally
knowing which irreducible $S^{(a,1^{p})}$ corresponds to can be
calculated directly, or also follows as a very special case of a theorem of
James \cite[Thm. A]{J2}.
\end{pf}

\newpage 

\section{Proof of the main result}

In this section we prove Theorem ~\reptypesuper.
One of the  main ideas will be to identify a subalgebra $eSe$ with
$k\mathcal{Q}/I$ for some quiver $\mathcal{Q}$. Then we can use
Gabriel's Theorem, since if the separated quiver is not a Dynkin or
Euclidean diagram then the algebra is of wild type.

\subsection{${\bf S(1|1, d):}$}

If $n=0$ or $m=0$ the superalgebra $S(m|n, d)$ is isomorphic to
the Schur algebra $S(m,d)$, so the smallest nontrivial case for us
is the superalgebra $S(1|1,d)$. This has been studied by Marko and
Zubkov:

\begin{prop}\cite[Prop. 2.1]{mz1}
\label{S11regrep} Let $S=S(1|1,d)$. If $p$ does not divide $d$
then $S$ is semisimple. If $p$ divides $d$ then $S$ has $d+1$
pairwise non-isomorphic one-dimensional modules, labelled by the
weights $0,1, \ldots , d$. The algebra is basic, connected, and
has left regular representation given by:
\begin{equation}
\label{leftregS11d} _SS=
\begin{array}{c}0\\1\\ \, \end{array}
\oplus
\begin{array}{c}1\\0\,\,\,\,\,\,2\\1\end{array}
\oplus
\begin{array}{c}2\\1\,\,\,\,\,\,3\\2\end{array}
\oplus \cdots \oplus
\begin{array}{c}r-1\\r-2\,\,\,\,\,\,\,\,\,\,\,\,\,\,\,\,r\\r-1\end{array}
\oplus
\begin{array}{c}r\\r-1\\ \, \end{array}
\end{equation}
\end{prop}

{}From here it is easy to obtain the quiver. In particular, for $m
\geq 1$ we define an algebra $\tilde{\mathcal{A}}_m$ to be the
algebra $k\mathcal{Q}/I$ where $\mathcal{Q}$ is the quiver with
$m+1$ vertices

\begin{picture}(428,50)
\put(80,20){\circle*{3}} \put(130,20){\circle*{3}}
\put(298,20){\circle*{3}}\put(348,20){\circle*{3}}

\put(85,22){\vector(1,0){40}}

\put(125,17){\vector(-1,0){40}} \put(135,22){\vector(1,0){40}}

\put(175,17){\vector(-1,0){40}}

\put(303,22){\vector(1,0){40}}

\put(343,17){\vector(-1,0){40}}

\put(253,22){\vector(1,0){40}} \put(293,17){\vector(-1,0){40}}

\put(190,20){\circle*{1}}

\put(210,20){\circle*{1}}\put(230,20){\circle*{1}}

\put(98,26){\footnotesize$\alpha_1$}
\put(98,7){\footnotesize$\beta_1$}
\put(148,26){\footnotesize$\alpha_2$}
\put(148,7){\footnotesize$\beta_2$}
\put(260,26){\footnotesize$\alpha_{m-1}$}
\put(260,7){\footnotesize$\beta_{m-1}$}
\put(316,26){\footnotesize$\alpha_{m}$}
\put(316,7){\footnotesize$\beta_{m}$}

\end{picture}

\noindent and $\mathcal{I}$ is the ideal generated by the
relations:
$$\alpha_i\alpha_{i+1}=0, \,\,\, \beta_{i+1}\beta_i=0
\,\,\,\alpha_1\beta_1=\beta_m\alpha_m=0 \mbox{ \rm and }
\beta_i\alpha_i=\alpha_{i+1}\beta_{i+1} \mbox{ \rm for } 1 \leq i
\leq m-1.$$

We remark that the quiver for $\tilde{\mathcal{A}}_m$ is the same
as the quiver for Erdmann's algebra (defined in \cite[3.1]{erd})
${\mathcal{A}}_{m+1}$ with one additional relation, namely that
$\beta_m\alpha_m=0$. As we will see in Section 5, this extra
relation is crucial in that it makes the algebra have infinite
global dimension.

\begin{thm}
\label{answerforS11} When $p$ divides $d$ the algebra $S(1|1,d)$
is isomorphic to the algebra $\tilde{\mathcal{A}}_d$ and has
finite representation type.
\end{thm}
\begin{pf}
The proof of the algebra isomorphism is essentially the same as
Proposition~3.2 in \cite{erd}, so we do not give it here. The algebra
$\tilde{\mathcal{A}}_m$ is a three nilpotent algebra, but
the only indecomposable modules of radical length three
are projective. Hence, the representation type of
$\tilde{\mathcal{A}}_m$ is the same as the representation type
of $\tilde{\mathcal{A}}_m$ modulo its radical squared (which
is two nilpotent). The separated quiver is a disjoint union of two
copies of the Dynkin diagram $A_{m+1}$, and hence by Gabriel's
Theorem the algebra is of finite type.
\end{pf}

\subsection{${\bf S(m|n,d),}$ ${\bf d<2p:}$ } We first
consider the case where $p \leq d <2p$. It turns out this
case is essentially the same as for the Schur algebra $S(2,d)$.
Recall that the summands of the signed permutation modules
$M^{(\alpha | \beta)}$ are called signed Young modules. We recall
the following theorem of Donkin.

\begin{thm}\cite[2.3 (6)]{Do}
\label{donkinsignedyoung} Let $m,n \geq d$. The isomorphism classes of indecomposible signed Young modules are labelled by the set
\[
\Lambda^{++}(m|n,d):=\{(\lambda | p\mu) \mid |\lambda|+p|\mu|=d\}
\] where $\lambda$ and $\mu$ are partitions.
\end{thm}  \noindent Following \cite{Do}, we write $Y^{(\lambda|p\mu)}$ for the indecomposible signed Young module labelled by $(\lambda|p\mu) \in \Lambda^{++}(m|n,d).$

We wish to show that for $d<2p$ the signed Young modules are
exactly the ordinary Young and twisted Young modules. To begin we
will show that the nonprojective signed Young modules are all
distinct.

Given a partition $\lambda \vdash d,$ let $\lambda'$ denote the transpose partition of $\lambda.$  For $\lambda \vdash d$ $p$-regular we let $m(\lambda)$ be the
Mullineux conjugate of $\lambda$, so $D^\lambda \otimes \sgn
\cong D^{m(\lambda)}$. It is well known that the Young module
$Y^{m(\lambda)'}$ is the projective cover of $D^\lambda$, a fact
we will use several times below.

\begin{prop}
\label{nonprojsigneddistinct} If the Young module $Y^\lambda$ is
isomorphic to a twisted Young module $Y^\mu \otimes \sgn$, then
$\lambda$ is $p$-restricted, in which case $\lambda =m(\mu')'$.
\end{prop}

\begin{pf} Applying the adjoint Schur functor to $Y^\lambda$ and $Y^\mu \otimes \sgn$
(see \cite[Thms 3.4.2, 3.8.2]{HN}) tells us that $P(\lambda) \cong
T(\mu')$. But the tilting modules are self dual, so $P(\lambda)$
is self-dual and hence has simple socle $L(\lambda)$. But
$V(\mu')$ is a submodule of $T(\mu')$ and the socle of $V(\mu')$
is $L(m(\mu')')$ so $\lambda = m(\mu')'$ is $p$-restricted.
\end{pf}

For $d<2p$ it is clear that the $p$-singular partitions are
exactly those of the form $(\mu_1+p, \mu_2, \ldots )$ for $\mu$ a
partition of $d-p$. Thus the number of signed Young modules for
$d$ is the total number of partitions of $d$ plus the number of
$p$-singular partitions, which by Prop. ~\nonprojsigneddistinct\
is precisely the number of Young and twisted Young modules. We
have shown the following result.

\begin{prop}
\label{nosignedyoung} For $p \leq d <2p$ the signed Young modules
are exactly the set of Young modules $Y^\lambda$ together with the
nonprojective twisted Young modules $Y^\mu \otimes \sgn$ where
$\mu$ is $p$-singular.
\end{prop}

We also remark that the non-projective Young modules are actually
irreducible and are in distinct blocks. Symmetric group blocks of
defect one are well understood and we have the following. Let
$\tau \vdash d-p$ and let $B$ be the block of $k\Sigma_d$ with
$p$-core $\tau$. Let

$$\lambda_1=(\tau_1+p, \tau_2, \ldots , \tau_s) \unrhd \lambda_2
\unrhd \cdots \unrhd \lambda_p=(\tau, 1^p)$$ be the $p$ partitions
of $d$ with $p$-core $\tau$.

Let $\mu = \tau'+(p)$. Then $\mu$ is $p$-singular and an easy
calculation shows $Y^\mu \otimes \sgn \cong D^{\lambda_{p-1}}$,
i.e. this is the nonprojective twisted Young module in the block.
So the signed Young modules in the block are:

\begin{equation}
\label{dlessthan2psignedyoung} Y^{\lambda_1}\cong D^{\lambda_1},
\,\,\,Y^{\lambda_2} \cong
\begin{array}{c}D^{\lambda_1}\\D^{\lambda_2}\\D^{\lambda_1}
\end{array}, \,\,\,
Y^{\lambda_i} \cong
\begin{array}{c}D^{\lambda_{i-1}}\\D^{\lambda_i}\,\,\,\,\,D^{\lambda_{i-2}}\\D^{\lambda_{i-1}}
\end{array},\,\,\,
Y^{\lambda_p} \cong
\begin{array}{c}D^{\lambda_{p-1}}\\D^{\lambda_{p-2}}\\D^{\lambda_{p-1}}
\end{array},\,\,\, Y^\mu \otimes \sgn \cong D^{\lambda_{p-1}},
\end{equation}
\noindent where $3 \leq i <p$. Thus the basic algebra for
$S(m|n,d)$ is a direct sum of two-sided ideals which are either
semisimple (corresponding to blocks of $k\Sigma_d$ of defect 0),
or which are isomorphic  to
\begin{equation}
\label{basic} \End_{k\Sigma_d}\left(
\bigoplus_{i=1}^pY^{\lambda_i}\oplus (Y^\mu \otimes \sgn)\right).
\end{equation}
The quiver and relations for the algebra in (\ref{basic}) is just
$\tilde{\mathcal{A}}_p$ which we saw earlier is of finite type. So
we can now prove the first half of Theorem ~\reptypesuper.

\begin{thm}
\label{answerfordlessthan2p} Let $S(m|n,d)$ be the
Schur superalgebra where $0\leq d <2p$.
\begin{itemize}
\item[(a)] $S(m|n,d)$ is semisimple if and only if one of the
following holds
\begin{itemize}
\item[(i)] $\operatorname{char} k=0$;
\item[(ii)] $d<p$;
\item[(iii)] $m=n=1$ and $p\nmid d$.
\end{itemize}
\item[(b)] $S(m|n,d)$ has finite representation type but is not
semisimple if and only if one of the following holds
\begin{itemize}
\item[(i)] $p\leq d < 2p$;
\item[(ii)] $m=n=1$ and $p\mid d$.
\end{itemize}
\end{itemize}
\end{thm}

\begin{pf} (a) If $d<p$ or $\operatorname{char } k=0$ then $k\Sigma_d$ is semisimple, hence so is the algebra
$S(m|n,d):=\End_{k\Sigma_d}(V^{\otimes d})$. Part (iii) follows by
Proposition~\S11regrep.

(b) The calculation above shows that $S(m|n,d)$ when
$\text{max}(m,n)\geq 2$ and $p\leq d <2p$ has finite
representation type (and is not semisimple). The same is true for
$S(1|1,d)$ when $p\mid d$ by Proposition~\S11regrep.
\end{pf}

\subsection{${\bf S(2|1,d), \,\,d\geq2p:}$}

In this section we begin to show that for $d \geq 2p$ the algebra
$S(2|1,d)$ has wild type. By Proposition ~\tensorbyonedim\ it will
suffice to prove this for $2p \leq d <3p$. Recall that:
$$S(2|1, d) \cong \End_{k\Sigma_d}\left( \bigoplus_{(\lambda_1, \lambda_2|\mu_1) \in
\Lambda(2|1,d)}M^{(\lambda_1, \lambda_2| \mu_1)}\right).$$ In each
case we will find some collection $\{Y_i\}$ of indecomposable
signed Young modules such that the algebra
$\End_{k\Sigma_d}(\oplus Y_i)$ is wild type. But this will imply
$S(2|1,d)$ is wild type by Proposition ~\descenttheorem.

To show that the algebra $\End_{k\Sigma_d}(\oplus Y_i)$ has wild
type we will show, in each case, that it is isomorphic to an
algebra $k\mathcal{Q} /I$ where $\mathcal{Q}$ is a quiver and $\mathcal{I}$
is an ideal contained in the ideal generated by paths of length at
least two. In each case the separated quiver of $\mathcal{Q}$ will
not be a union of Dynkin or extended Dynkin diagrams, which by
Gabriel's theorem means the algebra is of wild type. The main
difficulty will be to make an appropriate choice of the $\{ Y_i
\}$.

Our proof will split into three cases; each
following the same logic. First we handle the case
$2p \leq d <3p-2$. Then we handle $d=3p-2$ and $d=3p-1$
separately. That this happens is not surprising, as the blocks of
defect two are slightly different in the three cases, see for
instance Table 1 in \cite{Scopes}. Essentially everything is
understood for symmetric group blocks of defect two. The
decomposition numbers are all known by \cite{Richards}. The
structure of the projective indecomposable modules can be
determined from \cite{Scopes}. Finally Chuang and Tan \cite{ch}
have determined the structure of the nonprojective Young modules.
Thus we just need to judiciously select signed Young modules to
get a quiver of wild type. We will not give the details of
determining the module structures we need, they can be found in
the works cited above.

Recall that if $\lambda$ is a $p$-regular partition of $d,$ we write $m(\lambda)$ for the
Mullineux conjugate of $\lambda$, so $D^\lambda \otimes \sgn
\cong D^{m(\lambda)}$ and the Young module
$Y^{m(\lambda)'}$ is the projective cover of $D^\lambda.$

\subsection{${\bf S(2|1,d),}$ ${\bf 2p\leq d <3p-2:}$} Let $d=2p+t$
for $0 \leq t <p-2$. Using the Chuang-Tan notation we
have the following partitions in the principal block of
$\Sigma_d$:

$$\lambda^{(0)}=(2p+t),\,\,\,\,\,\lambda^{(1)}=(2p-1,t+1),\,\,\,\,\,\,
\lambda^{(2)}=(2p-2,t+1,1),\,\,\,\,\,\,\lambda^{(*)}=(p+t,p).$$

If we let $\sigma = (p+t,p-1,1)$ then $\sigma$ is $p$-restricted
and $Y^\sigma$ is the projective cover of $D^{\lambda^{(1)}}$.
Using $i$ to denote $D^{\lambda^{(i)}}$ we have the following
Loewy structures \cite[Thm. 2.4]{ch}:

\begin{equation}
\label{case1signedyoung}
 Y^{\lambda^{(0)}} \cong
0,\,\,\,\,\,Y^{\lambda^{(1)}} \cong
\begin{array}{c}0\\1\\0\end{array},\,\,\,\,\,Y^{\lambda^{(2)}}
\cong \begin{array}{c}1\\0\,\,\,\,\,\,\,2\\1\end{array},\,\,\,\,\,
Y^{\lambda^{(*)}}\cong
\begin{array}{c}1\\0\,\,\,\,\,\,\,*\\1\end{array},\end{equation}
\begin{equation}\nonumber Y^\sigma \cong
\begin{array}{c}1\\ *\,\,\,\,\,\,\,0\,\,\,\,\,\,\,2\\1\,\,\,\,\,\,\,\mu\,\,\,\,\,\,\,1\\
*\,\,\,\,\,\,\,0\,\,\,\,\,\,\,2\\1\end{array}.
\end{equation}
Since all the partitions above have third part zero or one and since by Proposition~\ab2trick
$Y^{(a,b,1)}$ is a summand of $M^{(a,b |1)}$, we know all the
modules in \eqref{case1signedyoung} are signed Young modules for $m=2, n=1.$
Consider the following quiver $\mathcal{L}$ given below:

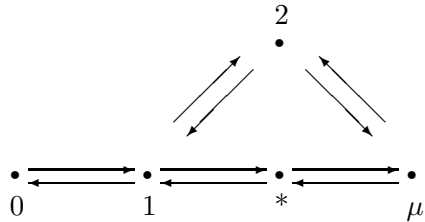
\begin{picture}(228,120)(-40,-20)

\put(80,20){\circle*{3}}
\put(78,5){0}

\put(130,20){\circle*{3}} \put(128,5){1}

\put(180,20){\circle*{3}}\put(178,5){*}

\put(230,20){\circle*{3}}\put(228,5){$\mu$}

\put(180,70){\circle*{3}}\put(178,77){2}

\put(85,22){\vector(1,0){40}}

\put(125,17){\vector(-1,0){40}} \put(135,22){\vector(1,0){40}}

\put(175,17){\vector(-1,0){40}}

\put(185,22){\vector(1,0){40}}

\put(225,17){\vector(-1,0){40}}


\put(140,40){\vector(1,1){25}}

\put(170,60){\vector(-1,-1){25}}

\put(220,40){\vector(-1,1){25}}

\put(190,60){\vector(1,-1){25}}

\end{picture}

\noindent Let $U$ be the direct sum of the five modules in
(\ref{case1signedyoung}). We have:

\begin{thm} \
\label{case1}
\begin{itemize}
\item[(a)] The algebra $\End_{k\Sigma_{2p+t}}(U)$ is isomorphic
to the algebra $k\mathcal{L}/I$ where I is in the ideal generated
by paths of length two. In particular it is of wild
representation type.
\item[(b)] The superalgebra $S(2|1, d)$ has wild
representation type for $2p \leq d <3p-2$.
\end{itemize}
\end{thm}
\begin{pf}
(a) The signed Young modules making up U have the same structure as
the modules $E_i$ for $1 \leq i \leq 5$ in Prop. 3.10 of
\cite{erd} and the determination of the quiver proceeds in exactly
the same way.

(b) This follows immediately from part (a) by
calculating the separated quiver of $\mathcal{L}$. The separated
quiver is two copies of:

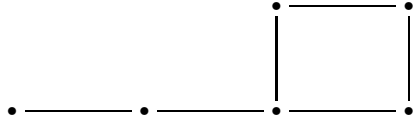
\begin{picture}(228,100)(-40,-20)

\put(80,20){\circle*{3}}

\put(130,20){\circle*{3}}

\put(180,20){\circle*{3}}

\put(230,20){\circle*{3}}
\put(230,60){\circle*{3}}
\put(180,60){\circle*{3}}

\put(180,24){\line(0,1){32}} \put(230,24){\line(0,1){32}}

\put(85,20){\line(1,0){40}}

\put(135,20){\line(1,0){40}}

\put(185,20){\line(1,0){40}}

\put(185,60){\line(1,0){40}}

\end{picture}

\noindent which is not a Dynkin or Euclidean diagram.
\end{pf}

\subsection{${\bf S(2|1,3p-2):}$} Let $d=3p-2$. Using the Chuang-Tan
notation we have the following
partitions in the principal block of $\Sigma_d$:

$$\lambda^{(0)}=(3p-2),\,\,\,\,\,\lambda^{(1)}=(2p-1,p-1),\,\,\,\,\,\,
\lambda^{(2)}=(2p-2,p-2,1^2),$$ $$\lambda^{(3)}=(2p-2,p-3,1^3),
\,\,\,\,\,\,\lambda^{(*)}=(2p-2,p),\,\,\,\,\,\,\,\delta=(2p-2,p-1,1).$$

\noindent It is also an easy calculation that:

$$m(\delta)' = (2p-3, p-1, 1^2)$$ so $Y^{(2p-3,p-1,1^2)}$ is
the projective cover of $D^\delta$. Let $\tau=(2p-3,p-1,1^2)$.
Using $i$ to denote $D^{\lambda^{(i)}}$ we have the following
structures for  Young modules \cite{ch}:

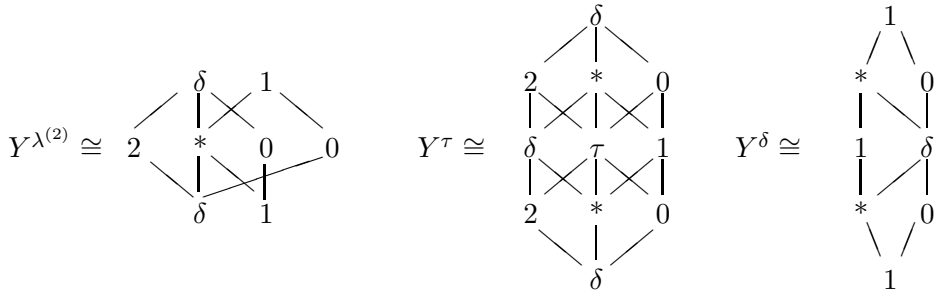
\begin{picture}(114,130)(-45,-30)

\put(-45,25){$Y^{\lambda^{(2)}}\cong$}

\put(3,36){\line(6,5){16}}

\put(31,36){\line(6,5){16}}

 \put(74,36){\line(-6,5){16}}
\put(48,36){\line(-6,5){16}}

\put(25,10){\line(-6,5){16}}

\put(48,12){\line(-6,5){16}}

\put(27,36){\line(0,1){13}}

\put(27,12){\line(0,1){13}}

\put(52,9){\line(0,1){13}}

\put(29,10){\line(3,1){45}}

\put(50,50){$1$}

\put(0,25){$2$}

\put(25,25){*}

\put(50,25){0}

\put(75,25){0}

\put(25,0){$\delta$}

\put(50,0){1}

\put(25,50){$\delta$}


\put(150,0){2} \put(175,0){*}\put(200,0){0} \put(150,25){$\delta$}
\put(175,25){$\tau$}\put(200,25){1}
\put(150,50){2}\put(175,50){*}\put(200,50){0}
\put(175,-25){$\delta$} \put(175,75){$\delta$}

\put(177.5,-14){\line(0,1){13}}

\put(177.5,11){\line(0,1){13}}

\put(177.5,36){\line(0,1){13}}

\put(177.5,61){\line(0,1){13}}

\put(152.5,11){\line(0,1){13}} \put(202.5,11){\line(0,1){13}}

\put(152.5,36){\line(0,1){13}} \put(202.5,36){\line(0,1){13}}

\put(158,36){\line(6,5){16}}

\put(158,11){\line(6,5){16}}

\put(158,61){\line(6,5){16}}

\put(183,36){\line(6,5){16}}

\put(183,11){\line(6,5){16}}

\put(183,-16){\line(6,5){16}}

\put(197,61){\line(-6,5){16}} \put(197,36){\line(-6,5){16}}
\put(197,11){\line(-6,5){16}}

\put(172,36){\line(-6,5){16}} \put(172,11){\line(-6,5){16}}
 \put(172,-16){\line(-6,5){16}}

\put(110,25){$Y^{\tau}\cong$}

\put(230,25){$Y^{\delta}\cong$}

\put(275,0){*}\put(300,0){0}

\put(275,25){1}\put(300,25){$\delta$}

\put(275,50){*}\put(300,50){0}

\put(286,75){1}

\put(286,-25){1}

\put(277.5,36){\line(0,1){13}}

\put(302.5,36){\line(0,1){13}}

\put(277.5,11){\line(0,1){13}}

\put(302.5,11){\line(0,1){13}}

\put(300,36){\line(-6,5){16}}

\put(300,24){\line(-6,-5){16}}

\put(280,60){\line(2,5){5}}

\put(280,0){\line(2,-5){5}}

\put(298,0){\line(-2,-5){5}}

\put(298,60){\line(-2,5){5}}
\end{picture}

\begin{equation}
\label{a=2youngmodules} Y^{\lambda^{(0)}} \cong 0, \,\,\,\,\,\,\,
Y^{\lambda^{(1)}} \cong
\begin{array}{c}0\\1\\0 \end{array},
\,\,\,\,\, \,\,
Y^{\lambda^{(3)}}\cong
\begin{array}{c}2\\\delta\,\,\,\,\,3\\2 \end{array}.
\end{equation}

All the Young modules in (\ref{a=2youngmodules}) are either for
two-part partitions or partitions of the form $(a,b,1^c)$ for
$c<p$, so by Proposition ~\ab2trick\ they are all signed Young
modules. Now consider the following quiver which we denote by
$\mathcal{Q}.$

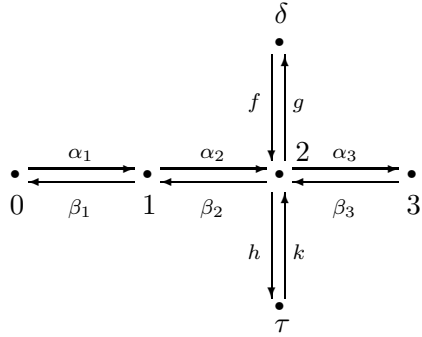
\begin{picture}(228,150)(-40,-50)
\put(80,20){\circle*{3}}
\put(78,5){0}

\put(130,20){\circle*{3}} \put(128,5){1}

\put(180,20){\circle*{3}}

\put(186,25){2}

\put(230,20){\circle*{3}}\put(228,5){3}

\put(180,70){\circle*{3}}\put(178,77){$\delta$}
\put(180,-30){\circle*{3}}\put(178,-40){$\tau$}

\put(182,25){\vector(0,1){40}} \put(177,65){\vector(0,-1){40}}
\put(182,-27){\vector(0,1){40}} \put(177,13){\vector(0,-1){40}}

\put(168,45){{\scriptsize $f$}}

\put(185,45){{\scriptsize $g$}}

\put(168,-12){{\scriptsize $h$}}

\put(185,-12){{\scriptsize $k$}}

 \put(85,22){\vector(1,0){40}}
 \put(100,26){{\scriptsize $\alpha_1$}}
  \put(150,26){{\scriptsize $\alpha_2$}}
    \put(200,26){{\scriptsize $\alpha_3$}}

 \put(100,5){{\scriptsize $\beta_1$}}
  \put(150,5){{\scriptsize $\beta_2$}}
    \put(200,5){{\scriptsize $\beta_3$}}

\put(125,17){\vector(-1,0){40}} \put(135,22){\vector(1,0){40}}

\put(175,17){\vector(-1,0){40}}

\put(185,22){\vector(1,0){40}}

\put(225,17){\vector(-1,0){40}}

\end{picture}

Let $T$ be the direct sum of the six signed Young modules in
(\ref{a=2youngmodules}). Then we have the following result.

\begin{thm}\
\label{case2}
\begin{itemize}
\item[(a)] The algebra $\Lambda=\End_{k\Sigma_{3p-2}}(T)$
is isomorphic to the algebra $k\mathcal{Q}/I$ where $\mathcal{I}$ is the ideal
generated the relations:
\begin{equation}
\label{relations} \begin{array}{c}
\alpha_1\alpha_2=\alpha_2\alpha_3=\alpha_2h
=\beta_3\beta_2=\beta_2\beta_1=\beta_3g=
f\alpha_3=\alpha_2h=0\\
hk=gf=\beta_1\alpha_2+\alpha_3\beta_3
\end{array}
\end{equation}
\item[(b)] The superalgebra $S(2|1, d)$ has wild
representation type for $d =3p-2$.
\end{itemize}
\end{thm}

\begin{pf}(a) Let $\{e_0, e_1, e_2, e_3, e_\delta, e_\tau\} \in
\Lambda$ be the canonical projections onto the Young module
summands of $T$. Choose nonzero maps $\alpha_i$ and $\beta_i$.
(The module structures leave no choice for the kernels and images
of these maps). The image and kernels of the maps $f, g, h, k$ can
also be determined by the structure of the modules given in
(\ref{a=2youngmodules}). In particular we have the following Loewy
structures for the kernels (from which the images are easily
determined):

$$\ker f \cong \begin {array}{c}*\,\,\,\,\,0\\1\end{array}, \,\,\,\,\ker k \cong
\begin{array}{c}\tau\\2\,\,\,\,\,*\,\,\,\,\,0\\\delta
\end{array},\,\,\,\,\,
\ker g \cong \begin{array}{c}2\,\,\,\,\,0\\\delta \end{array},
\,\,\,\,\, \ker h \cong \begin{array}{c}0\\1\end{array}
$$

\noindent Knowing the kernels and images of the maps is sufficient
to check that the zero relations in (\ref{relations}) are
satisfied. We need only show we can choose the maps so that
$hk=gf=\beta_1\alpha_2+\alpha_3\beta_3$.

Notice that the image of $\beta_2\alpha_2$ is the simple module 1
and the image of $\alpha_3\beta_3$ is the simple module $\delta$
while both $hk$ and $gh$ have image equal to the socle of
$Y^{\lambda^{(2)}}$. Since $\End (Y^{\lambda^{(2)}})$ is
three-dimensional, we can take as a basis $\{e_2, \alpha_3\beta_3,
\beta_2\alpha_2\}$. Write:
\begin{eqnarray*}hk&=&a\beta_2\alpha_2 + b \alpha_3\beta_3\\
gf&=&c\beta_2\alpha_2 + c \alpha_3\beta_3\end{eqnarray*} where the
images ensure $a,b,c,d$ are all nonzero.

Now just do the substitutions:
$$\beta_2'=a\beta_2, \,\, \alpha_2'=c\alpha_2,
\,\,\beta_3'=b\beta_3,\,\, \alpha_3'=d\alpha_4.$$

This will give the desired relation without effecting the earlier
zero relations. The maps $\{\alpha_i, \beta_i, h, k, f, g\}$
generate $\rad \Lambda$ and are independent modulo $I^2$, so there
is an epimorphism $k\mathcal{Q}/I$ onto $\Lambda$. Now one can
simply check the dimensions agree to obtain the isomorphism.

(b) This follows immediately from (a) by
calculating the separated quiver of $\mathcal{Q}$ and verifying
that it is not a Dynkin or Euclidean diagram.
\end{pf}

In the special case where $p=3$ (i.e. $S(2|1,7)$) the partition
$\lambda^{(3)}$ above does not exist. In this case to obtain the
same quiver just replace $Y^{\lambda^{(3)}}$ by the irreducible
module $D^\delta$, which is a signed Young module by Proposition
~\a0pcase. The quiver and relations are easily seen to be
identical.

\subsection{${\bf S(2|1,3p-1):}$}

Let $d=3p-1$. Using the Chuang-Tan notation we have the following
partitions in the principal block of $\Sigma_d$:

$$\lambda^{(0)}=(3p-1),\,\,\,\,\,\lambda^{(1)}=(2p-1,p-1,1),\,\,\,\,\,\,
\lambda^{(2)}=(2p-1,p-2,1^2),$$ $$\lambda^{(3)}=(2p-1,p-3,1^3),
\,\,\,\,\,\,\lambda^{(*)}=(2p-1,p),\,\,\,\,\,\,\,\rho=(2p-2,p,1),$$
\noindent where $\rho=m(\lambda^{(*)})'$, and hence $Y^\rho$ is
the projective cover of $D^{\lambda^{(*)}}$.

Using $i$ to denote $D^{\lambda^{(i)}}$ we have the following
structures for Young modules \cite{ch}:

\begin{equation}
\label{a=1youngmodules} Y^{\lambda^{(0)}} \cong
0,\,\,\,\,Y^{\lambda^{(*)}} \cong *,\,\,\,\,\,Y^{\lambda^{(1)}}
\cong
\begin{array}{c}0\,\,\,\,\,\, *\\1\\0\,\,\,\,\,\,*\end{array},
\,\,\,\,\,\,\,Y^{\lambda^{(2)}} \cong
\begin{array}{c}1\\0\,\,\,\,\,\,2\\1 \end{array},
\end{equation}

$$Y^{\lambda^{(3)}}\cong
\begin{array}{c}2
\\1\,\,\,\,\,\,3\\2 \end{array}, \,\,\,\,\ Y^\rho \cong
\begin{array}{c}
*\\
1\\
* \,\,\,\,\,\,0\,\,\,\,\,\, \rho\\
1\\
*
\end{array}.
$$
Now consider the following quiver called $\mathcal{H}$:

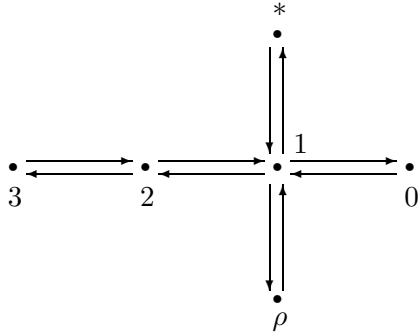
\begin{picture}(228,160)(-40,-60)
\put(80,20){\circle*{3}}
\put(78,5){3}

\put(130,20){\circle*{3}} \put(128,5){2}

\put(180,20){\circle*{3}}

\put(186,25){1}

\put(230,20){\circle*{3}}\put(228,5){0}

\put(180,70){\circle*{3}}\put(178,77){$*$}
\put(180,-30){\circle*{3}}\put(178,-40){$\rho$}

\put(182,25){\vector(0,1){40}} \put(177,65){\vector(0,-1){40}}
\put(182,-27){\vector(0,1){40}} \put(177,13){\vector(0,-1){40}}

 \put(85,22){\vector(1,0){40}}

\put(125,17){\vector(-1,0){40}} \put(135,22){\vector(1,0){40}}

\put(175,17){\vector(-1,0){40}}

\put(185,22){\vector(1,0){40}}

\put(225,17){\vector(-1,0){40}}

\end{picture}

Let $W$ be the direct sum of the six signed Young modules in
(\ref{a=1youngmodules}). Then we have the following result.

\begin{thm} \
\label{case3}
\begin{itemize}
\item[(a)] The algebra $\End_{k\Sigma_{3p-1}}(W)$ is
isomorphic to the algebra $k\mathcal{H}/I$ where $\mathcal{I}$ is in the ideal
generated by paths of length two. In particular it is of wild
representation type.
\item[(b)] The superalgebra $S(2|1, d)$ has wild
representation type for $d =3p-1$.
\end{itemize}
\end{thm}
\begin{pf}
(a) We leave the calculations here for the reader. They are similar to
the proof of Theorem~\caseB. (b) Once again the separated
quiver for $\mathcal{H}$ is not Dynkin or Euclidean so $S(2|1,3p-1)$
is of wild representation type.
\end{pf}

Just as in the $3p-2$ case, we remark that the partition
$\lambda^{(3)}$ does not exist when $p=3$ (i.e. $S(2|1,8)$).
Again, replacing $Y^{\lambda^{(3)}}$ with the simple module
$D^{(5,2,1)} \cong D^{\lambda^{(1)}}$ gives the same quiver.
$D^{(5,2,1)}$ is a signed Young module by Proposition ~\a0pcase.

\subsection{} We can now complete the proof of
Theorem~\reptypesuper\ by proving the following result.

\begin{thm} The algebra $S(m|n,d)$ has wild representation type
when $\operatorname{max}(m,n)\geq 2$ and $d\geq 2p$.
\end{thm}

\begin{pf} Since $S(m|n,d)\cong S(n|m,d)$ we may assume that
$m\geq n$. By assumption $m\geq 2$ and $n\geq 1$. The algebra
$S(2|1,d)$, $2p \leq d < 3p$ has wild representation type from
Theorems~\caseA (b), \caseB (b), \caseC (b). Hence, by
Proposition~\tensorbyonedim, $S(2|1,d)$ has wild type for $d\geq
2p$. One can now apply Corollary~\descentforschur(b) to conclude
that $S(m|n,d)$ has wild representation type for $d\geq 2p$.
\end{pf}

\section{Global Dimension}
\subsection{}
In \cite[Conj. 1]{mz1} it is conjectured that $S(m|n,d)$ is
quasi-hereditary whenever $d$ is coprime to $p$. In this section
we use our previous calculations to show this is far from true.  In fact, we will prove the following theorem.

\begin{thm}
\label{globaldimension} For $d \geq p \geq 5$ and $m,n \geq d$ the
superalgebra $S(m|n,d)$ has infinite global dimension.
\end{thm} This result suggests the following question.

\begin{quest} Is it true that $S(m|n,d)$ is either semisimple or has infinite global dimension? In particular, is it true that non-semisimple Schur superalgebras are
never quasi-hereditary?\label{conjecture}
\end{quest}

\subsection{}  Before proving Theorem~\globaldimension\ we first need to determine which block of $k\Sigma_{d}$ contains a given signed Young module.  Throughout the remainder of this section we assume $m,n \geq d.$

Recall the classification of irreducible $S(m|n,d)$-supermodules when $m,n \geq d$ given in \cite[2.3 (3)]{Do} (for the statement when $m,n,$ and $d$ are arbitrary, see \cite[Thm. 6.5]{bk}).

\begin{thm}
\label{donkinirreps} Let $m,n \geq d$. The irreducible supermodules of $S(m|n,d)$ are labeled by heighest weight by the set $\Lambda^{++}(m|n,d)=\{(\lambda |p\mu) \mid |\lambda|+p|\mu|=d\}$ where $\lambda$ and $\mu$ are partitions.
\end{thm} \noindent Write $L(\lambda|p\mu)$ for the irreducible $S(m|n,d)$-supermodule of highest weight $(\lambda|p\mu)$ and let $P(\lambda|p\mu)$ denote the projective cover of $L(\lambda|p\mu).$

Just as with the ordinary Schur algebra, one can define a Schur functor
\[
F: S(m|n,d)\text{-supermodules} \to k\Sigma_{d}\text{-modules}.
\] Like with the ordinary Schur functor, $F$ is exact and if $L$ is an irreducible $S(m|n,d)$-supermodule then $FL$ is either zero or an irreducible $k\Sigma_{d}$-module.  All irreducible $k\Sigma_{d}$-modules appear in this way.  In fact, by \cite[(5.13)]{bk} we have
\begin{equation}\label{schuronirreps}
FL(\lambda|p\mu) \cong \begin{cases} D^{m(\lambda')}, &\text{ if $|\mu|=0$ and $\lambda'$ is $p$-regular};\\
                                     0, & \text{ otherwise.}
\end{cases}
\end{equation}
 Furthermore, by \cite[p. 662]{Do} we have
\begin{equation}\label{schuronprojs}
FP(\lambda|p\mu) \cong Y^{(\lambda|p\mu)}.
\end{equation}
Thus to determine which block contains $Y^{(\lambda|p\mu)}$ it will suffice to obtain information about the composition factors of $P(\lambda|p\mu).$

Before proceeding, we require additional notation.  Let
\[
X=\left\{ \xi=(\xi_{1}, \dotsc , \xi_{m} | \xi_{m+1}, \dotsc , \xi_{m+n}) \mid  \xi_{1}, \dotsc , \xi_{m+n} \in \Z\right\}.
\]  We view $(\lambda|p\mu) \in \Lambda^{++}(m|n,d)$ as an element of $X$ in the natural way as $(\lambda, 0, \dotsc ,0 | p\mu, 0,\dotsc , 0).$
For $i=1, \dotsc ,m+n,$ let $\varepsilon_{i} \in X$ be the element which has a $1$ in the $i$th position and $0$'s elsewhere.  Let $\vartheta \in X$ be the element
\[
\vartheta=(-1, -2, \dotsc , -m| m-1, m-2, \dotsc , m-n+1,m-n).
\]  Note that our definition of $\vartheta$ differs from the one in \cite[(2.8)]{K} by a multiple of the $GL(m|n)$ analogue of the determinant representation.   This has no significant effect on our arguments.

Given $i=1, \dotsc , m+n$ we define $r_{i}: X \to \Z$ by
\[
r_{i}(\xi)=\begin{cases} \xi_{i}+\vartheta_{i}, & \text{ if $i=1, \dotsc , m$};\\
                         -(\xi_{i}+\vartheta_{i}), & \text{ if $i=m+1, \dotsc , m+n$}.

\end{cases}
\]  For $r \in \Zp$ and $\xi \in X,$ define
\begin{align*}
A_{r}(\xi)&=|\{i=1, \dotsc , m+n \mid r_{i}(\xi +\varepsilon_{i}) \equiv r \modp \}|, \\
B_{r}(\xi)&=|\{i=1, \dotsc , m+n \mid r_{i}(\xi) \equiv r \modp \}|.
\end{align*}

 For a partition $\lambda$ we let $\tilde{\lambda}$
denote its $p$-core.  Recall that the blocks of $k\Sigma_{d}$ are
parameterized by the $p$-cores of the partitions of $d$
\cite[Sect. 2.7]{JK}.  Also recall that the category of
$S(m|n,d)$-supermodules is equivalent to the category of
polynomial representations of degree $d$ for the supergroup
$GL(m|n).$ That is we can view the category of
$S(m|n,d)$-supermodules as a full subcategory of the category of
representations of $GL(m|n)$ and, in turn, as a full subcategory
of the category of supermodules for the superalgebra
$\operatorname{Dist}(GL(m|n))$ of distributions for $GL(m|n)$. See
\cite{bk} for a full discussion of these matters.
\begin{lem}
\label{signedblocks2} If $L(\lambda|p\mu)$ and $L(\nu|p\eta)$ have the same central character as $\operatorname{Dist}(GL(m|n))$-supermodules, then $\tilde{\lambda}=\tilde{\nu}.$
\end{lem}
\begin{pf}  If $L(\lambda|p\mu)$ and $L(\nu|p\eta)$ have the same central character, then by \cite[Lem. 3.3]{K} we have
\begin{equation}\label{E:ABequation}
A_{r}(\lambda|p\mu)-B_{r}(\lambda|p\mu)=A_{r}(\nu|p\eta)-B_{r}(\nu|p\eta)
\end{equation}
for all $r \in \Zp.$  However, note that if $i=m+1, \dotsc , m+n,$ then
\begin{align*}
r_{i}(\lambda|p\mu) &\equiv -\vartheta_{i} \modp \\
                    &\equiv r_{i}(\nu|p\eta) \modp,
\end{align*}
\noindent and
\begin{align*}
r_{i}((\lambda|p\mu)+\varepsilon_{i}) &\equiv -\vartheta_{i}-1 \modp \\
                    &\equiv r_{i}((\nu|p\eta)+\varepsilon_{i}) \modp.
\end{align*}  Consequently, we see that \eqref{E:ABequation} holds if and only if
\begin{equation}\label{E:ABprimeequation}
A'_{r}(\lambda|p\mu)-B'_{r}(\lambda|p\mu)=A'_{r}(\nu|p\eta)-B'_{r}(\nu|p\eta)
\end{equation}
for all $r \in \Zp,$ where for $\xi \in X$ we define
\begin{align*}
A'_{r}(\xi)&=|\{i=1, \dotsc , m \mid r_{i}(\xi +\varepsilon_{i}) \equiv r \modp \}|, \\
B'_{r}(\xi)&=|\{i=1, \dotsc , m \mid r_{i}(\xi) \equiv r \modp \}|.
\end{align*}

Given $(\lambda|p\mu)$ and $r \in \Zp,$ for brevity let us write $A'_{r}$ for $A'_{r}(\lambda|p\mu)$ and $B'_{r}$ for $B'_{r}(\lambda|p\mu).$  By definition, we have $A'_{r}=B'_{r-1}$ for all $r\in\Zp.$  We then have the following equations in the $B'_{r}$'s:
\begin{align*}
B'_{1}+B'_{2}+\dotsb +B'_{p}&=m \\
B'_{r-1}-B'_{r}&=:b_{r} \text{ (for $r=1, \dotsc , p$)}.
\end{align*}  However, a straightforward check verifies that the determinant for this system of linear equations is nonzero so one can solve this system for the $B'_{r}$'s.  That is, \eqref{E:ABprimeequation} holds if and only if
\begin{equation}\label{E:B'equation}
B'_{r}(\lambda|p\mu)=B'_{r}(\nu|p\eta)
\end{equation}
for all $r \in \Zp.$  At this point we assume without loss that $|\lambda|=|\nu|$ since $B'_{r}(\lambda,1^{p}|p\mu)=B'_{r}(\lambda|p\mu)$ and since the $p$-cores of $(\lambda, 1^{p})$ and $\lambda$ coincide.

Since for $i=1, \dotsc, m$ we have $r_{i}(\lambda|p\mu)=\lambda_{i}-i,$ it follows that  \eqref{E:B'equation} holds if and only if there is an element $\sigma \in \Sigma_{m}$ so that
\begin{equation}\label{E:nakequation}
\lambda_{i}-i \equiv \nu_{\sigma(i)} -\sigma(i) \modp
\end{equation}
for all $i=1, \dotsc, m.$
Finally, by the Nakayama Rule \cite[Thm. 5.1.1]{M} we have that \eqref{E:nakequation} holds if and only if $\tilde{\lambda}=\tilde{\nu}.$  This proves the desired result.
\end{pf}

\begin{cor}
\label{signedblocks}The signed Young module $Y^{(\lambda \mid
p\mu)}$ is in the $k\Sigma_{d}$ block with $p$-core $\tilde{\lambda}$.
\end{cor}
\begin{pf}  Let $L(\nu|p\eta)$ be a composition factor of $P(\lambda|p\mu).$  Then since both  $L(\lambda|p\mu)$ and $L(\nu|p\eta)$ are composition factors of an indecomposible $\operatorname{Dist}(GL(m|n))$-supermodule, they must have the same central characters (c.f. \cite[Sect. 2.8]{K}).  By the previous lemma, we have that $\tilde{\lambda}=\tilde{\nu}.$  The result then follows by exactness of the Schur functor along with \eqref{schuronirreps} and \eqref{schuronprojs}.
\end{pf}

\subsection{}  We can now prove Theorem~\globaldimension. Since $p \geq 5$ there must be a partition $\tau \vdash d-p$ which is a $p$-core
 by the work of Granville and Ono \cite{go}.   Consequently, $\Sigma_d$ has a $p$-block of defect
 one. Let $B$ denote this block.
\begin{lem} There are $p+1$ signed Young modules in $B$. They are
the $p$ distinct Young modules $Y^{(\lambda \mid \emptyset)} \cong
Y^\lambda$ with $\tilde{\lambda}=\tau$ and one non-projective twisted Young module $Y^{(\tau \mid p)} \cong Y^{(\tau' + p)} \otimes \sgn.$
\end{lem}
\begin{pf}
By Corollary~\signedblocks\ the signed Young modules which occur in $B$ are those of the form $Y^{(\lambda|p\mu)}$ with $\tilde{\lambda}=\tau$.  That is, precisely those listed.  Now, since twisted Young
modules are signed Young modules, we have $Y^{(\tau' + p)} \otimes
\sgn$ is a signed Young module in the block $B$, and it must be
$Y^{(\tau | p)}$.
\end{pf}

\noindent Now we have the direct sum decomposition

$$S(m|n,d)=\End_{k\Sigma_{d}}\left( V^{\otimes d} \right) \cong \End_{k\Sigma_{d}}\left( \bigoplus_{Y^{(\lambda \mid \mu)} \in
B}a_{(\lambda \mid \mu)}Y^{(\lambda \mid \mu)}\right) \oplus \End_{k\Sigma_{d}}(U)$$ where $a_{(\lambda \mid \mu)}$ denotes the multiplicity of $Y^{(\lambda|\mu)}$ in $V^{\otimes d}$ and $U$ has no summands in
the block $B$. That is, we have direct sum decomposition of $S(m|n,d)$ into graded two-sided ideals, one of 
which is Morita equivalent to
$$\End_{\Sigma_d}\left(\bigoplus_{Y^{(\lambda \mid \mu)} \in
B}Y^{(\lambda \mid \mu)}\right).$$ This is precisely the situation of
$(4.2.4)$ and the quiver is just $\tilde{\mathcal{A}}_{p}$. It is
easy to see this (super)algebra has infinite global dimension. The
projective resolutions of the simple modules are all periodic and
do not terminate. Since a two-sided graded ideal has infinite global
dimension, so does $S(m|n,d)$, proving Theorem ~\globaldimension.

\let\section=\oldsection

\end{document}